\documentclass[a4paper, 12pt, reqno]{amsart}
\usepackage{indentfirst}
\usepackage{graphicx}
\usepackage{amsmath, amsthm}
\usepackage{amssymb}
\usepackage{txfonts}
\usepackage{comment}

\theoremstyle{definition}
\newtheorem{Def}{Definition}[section]

\theoremstyle{theorem}
\newtheorem{Thm}[Def]{Theorem}

\newtheorem{Prop}[Def]{Proposition}

\newtheorem{ThmZ}{Theorem}

\theoremstyle{definition}
\newtheorem{Remark}{Remark}[section]

  \makeatletter
  \@addtoreset{equation}{section}
  \makeatother

\def\C{\mathbb{C}}
\def\sn{\mathrm{sn}}
\def\cn{\mathrm{cn}}
\def\dn{\mathrm{dn}}
\def\Re{\mathrm{Re}}
\def\Im{\mathrm{Im}}
\def\cos{\mathrm{cos}}
\def\sin{\mathrm{sin}}

\title[Variations of Bergman Kernels]{Variations of Bergman Kernels For Some Explicitly Given Families of planar Domains}
\author{Yanyan Wang}
\date{}

\thanks{2010 Mathematics Subject Classification. Primary 32A25; Secondary 32A17, 31A20.}

\address{
Graduate School of Mathematics\\
Nagoya University \\
Furo-cho, Chikusa-ku, Nagoya City, 464-8602, Japan }
\email{x08002k@math.nagoya-u.ac.jp}

\begin{document}

\maketitle {\bf Abstract:} { We study the parameter dependence of
the Bergman kernels on some planar domains depending on
complex parameter $\zeta$ in nontrivial ``pseudoconvex" ways. Smoothly bounded cases are studied at first: It
turns out that, in an example where the domains are
annuli, the Levi form for the logarithm of the Bergman kernels with respect to
$\zeta$ approaches to 0 as the point tends to the boundary of
the domain, and in another example where the domains are discs it approaches to $\infty$ as the point tends to the complement of a point in the boundary. Further, in contrast to this, in the cases where the boundary of the domains are not smooth, such as discs with slits, rectangles and half strips, completely
different phenomena are observed.
}
\section{Introduction and Results}
\par Bergman kernel and Bergman metric have been studied in detail in the case of bounded strongly pseudoconvex domains
with $C^{\infty}$ boundary. For such domains C. Fefferman
\cite{cf} (1974) found a remarkable asymptotic formula for the
Bergman kernel form. He used it to show that suitable geodesics of
the Bergman metric approach the boundary of the domain in a
``pseudotransverse'' manner, and therefore that biholomorphic
mappings of strongly pseudoconvex domains extend smoothly to the
boundaries. In 1978, using Fefferman's asymptotic formula, P.
Klembeck \cite{PK} showed, for a bounded strongly pseudoconvex
domain $\Omega$ with $C^{\infty}$ smooth boundary, that the
holomorphic sectional curvature of Bergman pseudometric near the
boundary of $\Omega$ approaches to the negative constant $-2/(n+1)$ which
 is the holomorphic sectional curvature of the
Bergman metric of the unit ball in $\C^n$. For some other results on the boundary behavior of the
Bergman kernel, refer to \cite{DH}, \cite{H}, \cite{M}.
\par In 2004, F. Maitani and H. Yamaguchi \cite{MY} brought a new viewpoint by studying the
variation of the Bergman metrics on the Riemann surfaces. Let us briefly recall it.
\par Let $B$ be a disk in the complex $\zeta$-plane, $\mathcal{D}$ be a domain in the product
space $B\times\mathbb{C}_{z}$, and let $\pi$ be the
first projection from $B\times\mathbb{C}_{z}$ to $B$ which is proper and smooth and
$D_\zeta=\pi^{-1}(\zeta)$ be a domain in $\C_z$.  Put $\partial
\mathcal{D}=\bigcup_{\zeta\in B}(\zeta,\partial D_\zeta)$.
\par Let $\zeta\in B$, $z \in D_\zeta$ and consider the
potential $\psi (\zeta,t,z)$ for $(D_\zeta,z)$ which is a complex valued harmonic function on $D_\zeta \backslash \{z\}$ vanishing on the boundary of $D_\zeta$, and decompose
$d_{t}\psi (\zeta,t,z)$ into
$$d_{t}\psi
(\zeta,t,z)=\mathrm{L}(\zeta,t,z)dt+\overline{\mathrm{K}(\zeta,t,z)}\overline{dt}$$
on $D_\zeta\backslash \{z\}$, where
$$\mathrm{L}(\zeta,t,z)=\frac{\partial \psi (\zeta,t,z)}{\partial
t},\mathrm{K}(\zeta,t,z)=\frac{\partial \overline{\psi
(\zeta,t,z)}}{\partial t}.$$
 Let $g(\zeta,t,z)$ be the Green function and $\lambda (\zeta,z)$  the Robin constant
for $(D_\zeta,z)$, then,
$$g(\zeta,t,z)=\log \frac{1}{|t-z|}+\lambda
(\zeta,z)+h(\zeta,t,z). $$ Here, $h(\zeta,t,z)$ is harmonic for
$t$ in a neighborhood of $z$ in $D_\zeta$ such that
$$h(\zeta,z,z)=0 \ \ \ \ \ \ \text{for}\ \ \zeta\in B.$$
Let $\varphi(\zeta,t)$ be a defining function of $ \partial\mathcal{D}$
in $B\times\mathbb{C}_{t}$ and define
$$k_{2}(\zeta,t)=(\frac{\partial ^{2}\varphi}{\partial \zeta\partial
\overline{\zeta}}|\frac{\partial \varphi}{\partial t}|^{2} -2\Re
\{\frac{\partial ^{2}\varphi}{\partial \zeta\partial
\overline{t}}\frac{\partial \varphi}{\partial \zeta}\frac{\partial
\varphi}{\partial \overline{t}}\}+|\frac{\partial
\varphi}{\partial \zeta}|^{2}\frac{\partial ^{2}\varphi}{\partial
t\partial \overline{t}})/|\frac{\partial \varphi}{\partial
t}|^{3}.$$
Then, the following variation formulae for the Bergman kernels $K_{D_\zeta}(z,z)$ on the diagonal can be obtained.

\begin{ThmZ}[\cite{MY}]
It holds for $(\zeta,z)\in \mathcal{D}$
\begin{eqnarray*}
\frac{\partial K_{D_\zeta}(z,z)}{\partial \zeta}&=&-\int\int_{D_\zeta}\frac{\partial \overline{\mathrm{K}(\zeta,t,z)}}{\partial\zeta}\mathrm{K}(\zeta,t,z)dxdy,\\
\frac{\partial^{2} K_{D_\zeta}(z,z)}{\partial \zeta\partial\overline{\zeta}}&=&
\frac{1}{4}\int_{\partial
{D_\zeta}}k_{2}(\zeta,t)(|\mathrm{L}(\zeta,t,z)|^{2}+|\mathrm{K}(\zeta,t,z)|^{2})ds_{t}\\
&+&\int\int_{D_\zeta}(|\frac{\partial
\mathrm{L}(\zeta,t,z)}{\partial \overline{\zeta}}|^{2}+|\frac{\partial
\mathrm{K}(\zeta,t,z)}{\partial \overline{\zeta}}|^{2})dxdy.
\end{eqnarray*}
\end{ThmZ}

Using Theorem 0.1, the following two properties of variation of
the Bergman kernels on the diagonal can be obtained.
\begin{ThmZ}[\cite{MY}] Let $\mathcal{D}$ be a pseudoconvex domain over
$B\times\mathbb{C}_{z}$ with smooth boundary, then $\log
K_{D_\zeta}(z,z)$ is plurisubharmonic on $\mathcal{D}$.
\end{ThmZ}
\begin{ThmZ}[\cite{MY}] Let $\mathcal{D}$
be a pseudoconvex domain over $B\times\mathbb{C}_{z}$
with smooth boundary. If, for each $\zeta\in B$, $\partial
\mathcal{D}$ has at least one strictly pseudoconvex point, then $\log K_{D_\zeta}(z,z)$ is a strictly plurisubharmonic
function on $\mathcal{D}$.
\end{ThmZ}
In 2006, B. Berndtsson \cite{B.B2} generalized the Theorem 0.2 to higher dimension and proved that,
\begin{ThmZ}[\cite{B.B2}] Let D be a pseudoconvex domain in $\C^k_\zeta\times \C^n_z$ and $\phi$ be a plurisubharmonic function on D. For each
$\zeta$ let $D_\zeta$ denote the $n$-dimensional slice
$D_\zeta:= \{z\in \C^n | (\zeta, z)\in D\}$
and by $\phi^\zeta$ the restriction of $\phi$ to $D_\zeta$. Let $K_{D_\zeta}(z,z)$ be the Bergman kernels of Bergman space $A_\zeta^2(D_\zeta, e^{-\phi^\zeta})$. Then, the function $\log K_{D_\zeta}(z,z)$ is plurisubharmonic
or identically equal to $-\infty$ on D.
\end{ThmZ}
\par These Theorems give us little information about the boundary behavior of the Bergman kernels on the diagonal, but, for any invariant metric on a domain $\Omega\subset \C^n$,
an important characteristic is its boundary behavior.
We now study the boundary behavior of the Levi form for the logarithm of Bergman kernels with respect to the parameter $\zeta$ by using complete explicit formulae of Bergman kernels on
certain families. This enables us to see precisely how
the Bergman kernels depends on the parameter $\zeta$.
\par Firstly, we consider two parameter domains with smooth boundaries, one is a family of discs, and the other is of annuli.
When the domains are annuli, we have the
following result.
\begin{Thm}
 Let $A_\zeta:=\{z\in\C\;|\;|\zeta|<|z|<1\},$
$\mathcal{A}:=\bigcup_{0<|\zeta|<1}\{\zeta\}\times A_\zeta$ and
$\partial \mathcal{A}:=
\bigcup_{0<|\zeta|<1}\{\zeta\}\times\partial A_\zeta$ then
$\partial^2 \log K_{A_\zeta}(z,z)/\partial \zeta\partial
\overline{\zeta}$ tends to $0$ with order 2 as $(\zeta,z)\in
\mathcal{A}$ tends to $\partial \mathcal{A}$ in a nontrivial way.
\end{Thm}\label{annuli}
 Here, the point $(\zeta,z)$ tends to the boundary in a nontrivial way means that $\zeta$
tends to a fixed point firstly, then $z$ tends to the boundary of the base domain $A_\zeta$. By parity of reasoning, we will repeat no more later.
\par The considered family of discs is
$$\mathcal{C}=\bigcup_{\zeta\in B}\{\zeta\}\times C_\zeta$$
where $B=\mathbb{D}_\zeta$,  $ C_\zeta=\{z\in\C_z|\, |z-e^{i\theta(\zeta)}|<1, \theta(0)=0,\, \Delta\theta(\zeta)=0\,\}$. Here, $\theta(\zeta)$ is a real-valued analytic function.
\par It is a well-known fact that a real hypersurface given by $|z-a(\zeta)|^2=e^{-\gamma(\zeta)}$ is Levi-flat if and only if
$$\gamma_{\zeta\overline{\zeta}}=2e^\gamma|a_{\overline{\zeta}}|^2, \ \ \ \ \ a_{\zeta\overline{\zeta}}+\gamma_{\zeta}a_{\overline{\zeta}}=0.$$ (For related results see B. Berndtsson \cite{B.B} and D. Barrett \cite{DB}.)
Therefore, $\partial \mathcal{C}$ is Levi-flat and we have the following theorem.
\begin{Thm}
The Levi form of $\log K_{C_\zeta}(z,z)$ with respect to $\zeta$ approaches to $\infty$ when $(\zeta,z)\in \mathcal{C}$ tends to $\partial\mathcal{C}\backslash\{(0,0)\}$ but depends on $\tan arg z$ when $(\zeta,z)\in \mathcal{C}$ tends to $(0,0)$ in a nontrivial way. When $(\zeta,z)$ tends to $(0,0)$ in a nontrivial way, if $\tan \arg z$ tends to $\infty$ then $\partial^2 \log K_{C_\zeta}(z,z)/\partial\zeta\partial\overline{\zeta}$ approaches to $\infty$, otherwise, $\partial^2 \log K_{C_\zeta}(z,z)/\partial\zeta\partial\overline{\zeta}$ approaches to a positive non-zero constant which depends on $\tan \arg z$.
\end{Thm}\label{circle}
\par Secondly, we investigate the boundary behavior of the Bergman kernels on particular pseudoconvex domain with non-smooth boundary. Since the variation formulae in Theorem 0.1 does not make
sense in the boundary when $\partial \mathcal{D}$ is not smooth, it is natural to
ask what will happen to the Bergman kernels on the diagonal in this case.
 We shall give an answer to this
question in a family of discs with slits.
\par The considered family of discs with slits are
$$D_{\zeta}=\{z\in\mathbb{D}_z\,|\,z\neq s\zeta,\, s\geq 1\,\},$$
where ${\zeta}\in B$ with
$B=\{\zeta\in\C|\,|\zeta-1|<\delta,|\zeta|<1\}$ and define
$\mathcal{D}=\bigcup_{\zeta\in B}\{\zeta\}\times D_{\zeta}.$
 Then the following result holds.
 \begin{Thm}The Levi form of $\log K_{D_\zeta}(z,z)$ with respect to $\zeta$ approaches to $\infty$ when $(\zeta,z)\in \mathcal{D}$ tends to $(1,1)\in\partial\mathcal{D}$ in a nontrivial way and approaches to $0$ when $(\zeta,z)$ tends to $(1,\pm i)\in \partial\mathcal{D}$ in a nontrivial way. Otherwise, $\partial^2 \log K_{C_\zeta}(z,z)/\partial\zeta\partial\overline{\zeta}$ approaches to a positive non-zero constant.
 \end{Thm}\label{discs}
\par On the other hand, the continuity of $K_{D_{\zeta}}(z,z)$
becomes already a delicate question in \cite{S 0}. So it also
might be interesting to see what will happen to the
Levi form of $\log K_{D_{\zeta}}(z,z)$ with respect to $\zeta$ as
$(\zeta,z)$ approach to a particular non-smooth boundary
point of a pseudoconvex domain with non-smooth boundary. Also, the explicit expression of Bergman kernel on a specific domain is a
 compelling problem. Although, there are three techniques for constructing the Bergman kernel, it's not easy to get the explicit
  expression of Bergman kernel on a specific domain. In this aspect, John P. D'Angelo \cite{JA}, D. Constales and R. S. Krau$\beta$har
   \cite{DC}, N. Suita \cite{NS} have done a lot of work.  We shall give the explicit expressions of Bergman kernels on rectangles by using the Schwarz-Christoffel transformation and Jacobi's elliptic functions, and on half strips by using the trigonometric functions.
\par The considered parameter rectangles
are
$$R_{\zeta}:=\{z=s+it\in \C_{z}\;|\;0<s<\Re\zeta,\,0<t<\Im\zeta
\}$$ where
 $\zeta\in B$ with $B:=\{\zeta\in\C||\zeta-(1+i)|<\delta\}$ and define
 $\mathcal{R}:=\bigcup_{\zeta\in B}\{\zeta\}\times R_{\zeta}.$ The angle of $\partial R_{\zeta} $ at each vertical point is $\pi/2$ and the following results hold.
\begin{Thm}The Bergman kernels of $R_{\zeta}$ on the diagonal
are
$$K_{R_{\zeta}}(z,z)=\frac{1}{\pi (\Im
\sn^2(u,k(\zeta)))^2}|\sn(u,k(\zeta))\cn(u,k(\zeta))\dn(u,k(\zeta))\frac{K(k(\zeta))}{\Re\zeta}|^2$$
where $u=K(k(\zeta))z/\Re\zeta$ and
$\sn(u,k),\cn(u,k),\dn(u,k)$ are the Jacobi's elliptic funtions of
the first kind, $K(k)$ is the complete elliptic integral
of the first kind. $k(\zeta)$ is a real valued analytic function
with respect to $\zeta$. Let $\zeta=1+i+\varepsilon$ then,
$$k(\zeta)=k_0+2\Re((a+ib)\varepsilon)+2\Re((c+id)\varepsilon^2)+2e|\varepsilon|^2+\cdots,$$
where $k_0=1/\sqrt{2},a=b=-2c=K/\left(4\sqrt{2}(2E-K)\right)$,
$d=e=-\sqrt{2}a^2$. Here, $K$ is the value of the complete elliptic integral
of the first kind at the point $k=1/\sqrt{2}$, and $E$ is the
value of the complete elliptic integral of the second kind at the
point $k=1/\sqrt{2}$.
\end{Thm}\label{rectangles}
\begin{Thm} For Bergman kernels $K_{R_{\zeta}}(z,z)$
where $(\zeta,z)\in\mathcal{R}$, it holds that
\begin{eqnarray*}
 \lim_{z\rightarrow 0}\lim_{{\zeta}\rightarrow
1+i}\frac{\partial^2\log K_{R_{\zeta}}(z,z)}{\partial {\zeta}
\partial\overline{\zeta}}=0, \ \ \ \ \ \ \ \  \ \ \ \lim_{z\rightarrow 1}\lim_{{\zeta}\rightarrow 1+i}\frac{\partial^2\log K_{R_{\zeta}}(z,z)}{\partial
\zeta
\partial\overline{\zeta}}=\frac{3}{2},\\
\lim_{z\rightarrow i}\lim_{{\zeta}\rightarrow 1+i}\frac{\partial^2\log K_{R_{\zeta}}(z,z)}{\partial
\zeta
\partial\overline{\zeta}}=16a^2,\ \ \ \ \ \lim_{z\rightarrow 1+i}\lim_{{\zeta}\rightarrow 1+i}\frac{\partial^2\log K_{R_{\zeta}}(z,z)}{\partial
\zeta
\partial\overline{\zeta}}=\infty.
\end{eqnarray*}
\end{Thm}\label{Rectangles}
In the case of half strips, we get the following Theorem.
\begin{Thm}Let $S_\zeta:=\{z\in\C_z|0<\Re
z<\Re\zeta,\,\Im z>0\},$ with $B=\{\zeta\in\C||\zeta-1|<\delta\}$
and $\mathcal{S}:=\bigcup_{\zeta\in B}\{\zeta\}\times S_{\zeta}.$
Then,
$$\lim_{z\rightarrow 0}\left(\lim_{\zeta\rightarrow
1}\frac{\partial^2\log K_{S_\zeta}(z,z)}{\partial\zeta\partial\overline{\zeta}}\right)=0,\ \ \lim_{z\rightarrow 1}\left(\lim_{\zeta\rightarrow
1}\frac{\partial^2\log K_{S_\zeta}(z,z)}{\partial\zeta\partial\overline{\zeta}}\right)=\infty.
$$
But,
$$\lim_{\Im z\rightarrow \infty}\left(\lim_{\zeta\rightarrow
1}\frac{\partial^2\log K_{S_\zeta}(z,z)}{\partial\zeta\partial\overline{\zeta}}\right)$$
depends on the choice of  $\Re z$.
\end{Thm}\label{strips}
\begin{Remark}
Theorem 1.5 and 1.6 indicate that the Levi form of $\log K_{D_\zeta}(z,z)$ with respect to $\zeta$ is independent of the choice of the parameter $\zeta$ at the point $(\zeta,0)$, but it is sensitively dependent on the parameter
$\zeta$ on the other singular points.
\end{Remark}

\section{Preliminaries}
 \def\theequation{\thesubsection.\arabic{equation}}

\par We briefly present here certain results underlying the proofs
of Theorems. This exposition is adapted to our special cases.
\subsection{Bergman Kernel}
 The Bergman kernel of a domain $\Omega\subset\C^n$ is a reproducing
kernel for the Hilbert space of all square integrable holomorphic
functions on $\Omega$. In what follows, let $\Omega$ be a bounded
domain in $ \mathbb{C}^{n}$, let $A^{2}(\Omega)$ be the space of
square integrable holomorphic functions on $\Omega$.
And let $\{\phi_{j}(z) \}_{j=0}^{\infty}$ be a complete
orthonormal basis for $A^{2}(\Omega)$. Then the Bergman kernel
$K_{\Omega}(z,w)$ is identified with the following series:
$$K_{\Omega}(z,w)=\sum_{j=0}^{\infty}{\phi_{j}(z) }\overline{{\phi_{j}(w)
}},$$ which is independent of the choice of orthonormal basis. For
$z = w,$ one has $K_{\Omega}(z,z)>0$.
\par The Bergman kernel satisfies the following transformation formula.
\begin{Prop}Let $f:\;\Omega\longrightarrow D $ be a
biholomorphic mapping between $\Omega$ and $D$. Then,
$$K_{\Omega}(z,w)=K_{D}(f(z),f(w))\det f'(z)\overline{\det f'(w)}.$$
\end{Prop}
 By Cauchy's estimate it is easy to see that $K_{\Omega}(z,w)$ is a $C^{\infty}$¡¡function
on $\Omega \times\Omega$ and on the diagonal, it can be
represented as
 $$K_{\Omega}(z,z)=\sup \{|f(z)|^{2}\mid\,f\in A^{2}(\Omega),||f(z)||_{A^{2}(\Omega)}=1 \} \;\; \text{for} \;\;\forall\;z\in
 \Omega.
 $$
 \subsection{Schwarz-Christoffel Transformation}
 Riemann's mapping theorem states that any non-empty simply connected domain which is neither the whole nor the extended $z$ plane can be mapped with a univalent transformation onto the unit disk or onto the upper half of the complex $w$ plane. Unfortunately, there is no general constructive approach for finding the univalent transformation. Nevertheless, as we will see, there are many particular domains, such as the interior of a polygon, for which the univalent function can be constructed explicitly.
\par Let $\Gamma$ be a piecewise linear boundary of a polygon in
 the $w$-plane and let the interior angles at successive vertices
 be $\alpha_{1}\pi,\alpha_{2}\pi,\cdots,\alpha_{n}\pi$. The
 transformation defined by the equation
$$w=F(z)=C\int_{0}^{z}(\xi-a_{1})^{\alpha_{1}-1}(\xi-a_{2})^{\alpha_{2}-1}\cdots(\xi-a_{n})^{\alpha_{n}-1}d\xi+C',\eqno (2.2.1)$$
where $C,\,C'$ are complex numbers and $a_{1},a_{2},\cdots,a_{n}$
are real numbers, maps $\Gamma$ into the real axis of the complex
$z$ plane and the interior of the polygon to the upper half of the
$z$ plane. The vertices of the polygon $A_{1},A_{2},\cdots A_{n}$
are mapped to the points $a_{1},a_{2},\cdots,a_{n}$. This map is
an analytic one-to-one conformal transformation between the upper
half of the $z$ plane and the interior of the polygon.

\begin{Remark} Actually, for any univalent transformation, the
correspondence of three points on the boundaries of two simply
connected domains can be prescribed arbitrarily. In particular,
any of the three vertices of the polygon can be associated with
any three points on the real axis.
\end{Remark}
\par An interesting application of the Schwarz-Christoffel construction is the mapping of a
rectangle. Despite the fact that it is a simple closed polygon,
the function defined by Schwarz-Christoffel transformation is not
elementary. In the case of a rectangle, we find that the mapping
functions involve elliptic integrals and elliptic functions.
\subsection{Elliptic Integrals}
The incomplete elliptic integral of the first kind is
$$u=F(w,k):=\int_0^w\frac{dt}{\sqrt{(1-t^2)(1-k^2t^2)}}.$$
The parameter $k$ is called the modulus of the elliptic integral (for more information about $k$, refer to \cite{s 1}).
$$K(k):=F(1,k)=\int_0^1
\frac{dt}{\sqrt{(1-t^2)(1-k^2t^2)}}=\int_0^{\frac{\pi}{2}}\frac{d\theta}{\sqrt{1-k^2\sin^2\theta}},$$
which is referred to as the complete elliptic integral.
$K(k)$ has a power series expansion
$$K(k)=\frac{\pi}{2}\left(1+\sum_{n=1}^\infty\frac{((2n-1)!!)^2}{((2n)!!)^2}k^{2n}\right).\;\;\;\;\;\;\;\;\;\;\;\;\;\;\;\;\;\;$$
 The special values of $K(k)$ are
$K(0)=\pi/2,\;K(1)=\infty,\;K(2^{-\frac{1}{2}})=\pi^{-\frac{1}{4}}\Gamma(1/4)^2/4$.
\par The incomplete elliptic integral of the second kind is defined
by
$$E(u,k)=\int_0^u\dn^2(w,k)\;dw.\;\;\;\;\;\;\;\;\;\;\;\;\;\;\;\;\;\;\;\;\;\;$$
 Its power series expansion with respect to
 $u$ in the neighborhood of $u=0$ is
 $$E(u,k)=u-\frac{k^2}{6}u^3+\frac{4k^2-3k^4}{5!}u^5+O(u^7).$$
The complete elliptic integral of the second kind is defined by
$$E(k)=E(K(k),k)=\int_0^1\sqrt{\frac{1-k^2t^2}{1-t^2}}\;dt=\int_0^{\frac{\pi}{2}}\sqrt{1-k^2\textmd{sin}^2\theta}\;d\theta.$$
The derivatives of $K(k)$ and $E(k)$ with respect to $k$ are
$$\frac{\partial}{\partial
k}K(k)=\frac{1}{k(1-k^2)}\left(E(k)-(1-k^2)K(k)\right),\;\;\;\;\;\;\;\frac{\partial}{\partial k}E(k)=\frac{1}{k}\left(E(k)-K(k)\right).\;\;\;\;\;\;\;\;\;\;\;\;\;\;\;\;\;\;$$
The power series expansion of $E(k)$ is
$$E(k)=\frac{\pi}{2}\left(1-\sum_{n=1}^\infty\frac{((2n-1)!!)^2}{((2n)!!)^2}\frac{k^{2n}}{2n-1}\right).\;\;\;\;\;\;\;\;$$
The special values of $E(k)$ are $E(0)=\pi/2,\;E(1)=1$, and
$E(2^{-\frac{1}{2}})=\pi^{\frac{3}{2}}\Gamma(1/4)^{-2}+\pi^{-\frac{1}{2}}\Gamma(1/4)^2/8$.
\subsection{Elliptic Functions} In this section we introduce some kinds of elliptic functions.
\subsubsection{Weierstrass's elliptic functions} For
$\omega_1$, $\omega_2\in\C$, linearly independent over
$\mathbb{R}$, let
$\Lambda:=\{2j\omega_1+2k\omega_2|(j,k)\in\mathbb{Z}^2\}$ be the
lattice in $\C$. We define the Weierstrass elliptic function
$\mathcal{P}$ by
$$\mathcal{P}(z):=\mathcal{P}(z;\omega_1,\omega_2)=\frac{1}{z^2}+\sum_{0\neq\omega\in\Lambda}\left(\frac{1}{(z-\omega)^2}-\frac{1}{\omega^2}\right).$$
The periodicity of $\mathcal{P}$ is
$$\mathcal{P}(z+2\omega_1)=\mathcal{P}(z)=\mathcal{P}(z+2\omega_2).$$
Weierstrass's zeta function is determined by $\zeta'(z)=-\mathcal{P}(z)$, then
$$\zeta(z)=\frac{1}{z}-\sum_{0\neq\omega\in\Lambda}\left(\frac{1}{z-\omega}+\frac{z}{\omega^2}+\frac{1}{\omega}\right).$$
\subsubsection{Jacobi's elliptic functions of the first kind}
 The Jacobi's elliptic functions of the first kind are
$\sn(u,k)$, $\cn(u,k)$ and $\dn(u,k)$. A standard ``normalized"
definition of $\sn(u,k)$ is
$$w=F^{-1}(u,k)=\sn(u,k).$$
Then the functions $\cn(u,k)$ and $\dn(u,k)$ can now be defined,
respectively, by
\begin{equation}
\cn^2(u,k)=1-\sn^2(u,k),\;\;\;\;\;\; \dn^2(u,k)=1-k^2\sn^2(u,k).
\end{equation}
$\sn u$, $\cn u$, $\dn u$ are double periods functions and satisfy
\begin{eqnarray}
\sn (u+4K(k))&=&\sn(u+2iK'(k))=\sn u,\\
\cn(u+4K(k))&=&\cn(u+2K(k)+2iK'(k))=\cn u,\\
\dn(u+2K(k))&=&\dn(u+4iK'(k))=\dn u,
\end{eqnarray}
where $K'(k)=K(k')$, here $k'=\sqrt{1-k^2}$ is the complementary modulus.
We will use the following special values of $\sn u$, $\cn u$, $\dn u$ in the proof of Theorem 1.5.
\begin{eqnarray}
\sn K(k)&=&1,\ \ \ \ \sn iK'(k)=\infty,\ \ \ \ \sn (K(k)+iK'(k))=1/k,\\
\cn K(k)&=&0,\ \ \ \ \cn iK'(k)=\infty,\ \ \ \ \cn (K(k)+iK'(k))=k'/ik,\\
\dn K(k)&=&k',\ \ \ \ \cn iK'(k)=\infty,\ \ \ \ \dn (K(k)+iK'(k))=0.
\end{eqnarray}
 Just as the trigonometric functions have simple algebraic addition theorems,
$\sn (u,k)$, $\cn (u,k)$ and $\dn(u,k)$ also have addition theorem as below.
\begin{eqnarray}
\sn(u+v)&=&\frac{\sn u\cn v\dn v+\sn v\cn u\dn u}{1-k^2\sn^2u\sn^2v},\\
\cn(u+v)&=&\frac{\cn u\cn v-\sn u\sn v\dn u\dn v}{1-k^2\sn^2u\sn^2v},\\
\dn(u+v)&=&\frac{\dn u\dn v-k^2\sn u\sn v\cn u\cn v}{1-k^2\sn^2u\sn^2v}.
\end{eqnarray}
Putting $u=v$, the following duplication formulae result,
\begin{eqnarray}
\sn 2u&=&\frac{2\sn u\cn u\dn u}{1-k^2\sn^4 u},\\
\cn 2u&=&\frac{2\cn^2 u}{1-k^2\sn^4 u}-1=1-\frac{1\sn^2u\dn^2u}{1-k^2\sn^4 u},\\
\dn 2u&=&\frac{2\dn^2 u}{1-k^2\sn^4 u}-1=1-\frac{2k^2\sn^2u\cn^2u}{1-k^2\sn^4 u}.
\end{eqnarray}
From these results we can verify that the formulae listed below hold,
\begin{equation}
\sn u=\sqrt{\frac{1-\cn 2u}{1+\dn 2u}},\ \ \cn u=\sqrt{\frac{\cn 2u+\dn 2u}{1+\dn 2u}},\ \ \dn u=\sqrt{\frac{\cn 2u+\dn 2u}{1+\cn 2u}}.
\end{equation}
 In the neighborhood of $u=0$, the power series expansions for the Jacobi's elliptic functions of the first kind with respect to $u$ are given below:
\begin{eqnarray*}
\sn(u,k)&=&u-\frac{1}{3!}(1+k^2)u^3+\frac{1}{5!}(1+14k^2+k^4)u^5+O(u^7),\\
\cn(u,k)&=&1-\frac{1}{2!}u^2+\frac{1}{4!}(1+4k^2)u^4+O(u^5),\\
\dn(u,k)&=&1-\frac{1}{2!}k^2u^2+\frac{1}{4!}(4k^2+k^4)u^4+O(u^6).
\end{eqnarray*}
We refer to \cite{s 1} for more information on elliptic integrals and elliptic functions.
\section{Proofs of Theorems and Remarks }We next accomplish the proofs of Theorems that are given in
the introduction.
\subsection{Case of The Family of Annuli}
\makeatother
\makeatletter
\renewenvironment{proof}[1][\proofname]{\par
  \pushQED{\qed}%
  \normalfont \topsep6\p@\@plus6\p@\relax
  \trivlist
  \item[\hskip\labelsep
        \textrm
    #1\@addpunct{.}]\ignorespaces
}{%
  \popQED\endtrivlist\@endpefalse
}
\makeatother

\begin{proof}[Proof of Theorem 1.1]
It is well known from \cite{NS} that
$$
 K_{A_\zeta}(z,z)=\frac{\mathcal{P}(u)+c(\omega_1)}{\pi |z|^2},
\eqno(3.2.1)$$
where $\omega_1=-\log|\zeta|$, $u=-2\log |z|\in (0,2\omega_1)$,
$c(\omega_1)=\zeta(\omega_1)/\omega_1$ and
 $\mathcal{P}(u)$ is the Weierstrass elliptic
function with periods 2$\omega_1$, 2$\omega_2=2\pi i$, $\zeta(u)$
is the Weierstrass's zeta function.
\par In addition, from (3.2.1) we can get that
$$\frac{\partial^2 }{\partial \zeta\partial \overline{\zeta}}\log K_{A_\zeta}(z,z)
=e^{2\omega_1}\frac{\left(2\mathcal{P}(u)-\mathcal{P}(\omega_1)+c\right)(\mathcal{P}(\omega_1)+c)}
{4\omega_1^{2}(\mathcal{P}(u)+c)^2}.$$
$\mathcal{P}(0)=\infty$ and $\mathcal{P}(u)$ decreases in
$(0,\omega_1)$ can be easily checked . Also, we know that
$\mathcal{P}(2\omega_1-u)=\mathcal{P}(u)$ and
$\omega_1^2\mathcal{P}(\omega_1)=\pi^2/6$. So
$\mathcal{P}(u)>0$ in $(0,2\omega_1)$ can be obtained.
  And note that in the neighborhood of
$u=0$, $\mathcal{P}(u)$ has the second order pole, that is,
$$\mathcal{P}(u)=u^{-2}(1+O(u^2)).$$
 Then,
$$2\mathcal{P}(u)-\mathcal{P}(\omega_1)+c=2u^{-2}(1+O(u^2)),$$
$$(\mathcal{P}(u)+c)^2=u^{-4}(1+O(u^2)),$$
in the neighborhood of $u=0$. If $|z|\rightarrow 1$ then $u
\rightarrow 0$.
Hence,
$$\lim_{|z|\rightarrow 1}\frac{\partial^2 \log
K_{A_\zeta}(z,z)}{\partial \zeta\partial \overline{\zeta}}=0$$
with order 2. Using the periodicity of $\mathcal{P}(u)$,
similarly, we also can deduce that
$$\lim_{|z|\rightarrow |\zeta|}\frac{\partial^2 \log
K_{A_\zeta}(z,z)}{\partial \zeta\partial \overline{\zeta}}=0$$
with order 2.
\end{proof}
\begin{Remark}
 The proof of Theorem 1.1 implies that although the Levi form of $\log K_{D_\zeta}(\zeta,z)$ with respect to $\zeta$ approaches to 0 when $(\zeta,z)$ tends to the boundary of the domain,  $\log K_{D_\zeta}(\zeta,z)$ is a strictly plurisubharmonic function on $\mathcal{A}$.
\end{Remark}
\subsection{Case of The Family of Discs}
\begin{proof}[Proof of Theorem 1.2]
Using the Proposition 2.1, we can get the Bergman kernels
$$K_{C_\zeta}(z,z)=\frac{1}{\pi\left(1-|z+e^{i\theta(\zeta)}|^2\right)^2}.$$
 So, $$\frac{\partial^2\log K_{C_\zeta}(z,z)}{\partial \zeta\partial \overline{\zeta}}=4|\theta_{\zeta}|^2\frac{|z|^2\left(\Re (ze^{-i\theta(\zeta)})+2\right)}{\left(|z|^2+2\Re (ze^{-i\theta(\zeta)})\right)^2}.$$
 Moreover, the condition $\theta(0)=0$ induces that
 $$
  \lim_{\zeta\rightarrow 0}\frac{\partial^2\log K_{C_\zeta}(z,z)}{\partial \zeta\partial \overline{\zeta}}=4|\theta_{\zeta}|^2\frac{|z|^2(\Re z+2)}{\left(1-|z+1|^2\right)^2}.
\eqno(3.1.1)$$
  Then, from (3.1.1), if $(\zeta,z)\in \mathcal{C}\,$ tends to $\partial\mathcal{C}\backslash\{(0,0), (0,-2)\}$ then $\partial^2 \log K_{C_\zeta}(z,z)/\partial\zeta\partial\overline{\zeta}$ tends to $\infty$. Also, from (3.1.1) we have
  \begin{equation*}
 \lim_{z\rightarrow -2} \lim_{\zeta\rightarrow 0}\frac{\partial^2\log K_{C_\zeta}(z,z)}{\partial \zeta\partial \overline{\zeta}}=\lim_{z\rightarrow -2}\lim_{\zeta\rightarrow 0}4|\theta_{\zeta}|^2\frac{1}{\left(2+\Re z\right)^2}=\infty.
  \end{equation*}
 That is, if $(\zeta,z)\in \mathcal{C}$ tends to $(0,-2)\in\partial\mathcal{C}$ then $\partial^2 \log K_{C_\zeta}(z,z)/\partial\zeta\partial\overline{\zeta}$ tends to $\infty$ with order $1$.
 \par Finally, we consider the boundary point $(0,0)$.
   Let $z=x+iy$, then
$$\lim_{\zeta\rightarrow 0}\frac{\partial^2\log K_{C_\zeta}(z,z)}{\partial \zeta\partial \overline{\zeta}}=4|\theta_{\zeta}|^2\frac{(x+2)}{\left(x^2+y^2+4x+\frac{4}{\left(1+(y/x)^2\right)}\right)}.$$
Therefore,
 when $(\zeta,z)$ tends to $(0,0)$, the Levi form of $\log K_{C_\zeta}(z,z)$ with respect to $\zeta$ is dependent on $\tan \arg z$ and if $\tan \arg z$ tends to $\infty$ the Levi form of $\log K_{C_\zeta}(z,z)$ with respect to $\zeta$ tends to $\infty$, otherwise, the Levi form of $\log K_{C_\zeta}(z,z)$ with respect to $\zeta$ tends to a positive non-zero constant which depends on $\tan \arg z$.
\end{proof}
\subsection{Case of The Family of Discs with Slits}
  We will now proceed with the proof of Theorem 1.3.
\begin{proof}[Proof of Theorem 1.3]
\par Define that
$$z=E_{\zeta}(w):=e^{-i\theta}K^{-1}(e^{-t}K(w)),$$
where $K(w)=w/(1+w)^{2}$ is the Koebe function,
${\zeta}=e^{-i\theta}K^{-1}(e^{-t}/4)$ with $t>0,\,
\theta\in[0,2\pi)$, this is inspired by \cite{s 3}.
 For each $\zeta$,  $E_{\zeta}(w)$ maps the unit disc $\mathbb{D}$ in the $w$ plane to $D_{\zeta}$ which is the unit disc in the complex $z$
plane minus a segment $[\zeta,\,e^{-i\theta})$.
  Then the inverse mapping of $E_{\zeta}(w)$
is
$$w=E_{\zeta}^{-1}(z)=K^{-1}(e^{t}K(e^{i\theta}z)).$$
 We now can get
the power series expansion of $E_{\zeta}^{-1}(z)$ with respect to
the parameter $\zeta$ in a neighborhood of $\zeta=1$, that is, in
a neighborhood of $(t,\theta)=(0,0)$ which is
\begin{eqnarray*}
 E_{\zeta}^{-1}(z)&=&K^{-1}(e^{t}K(e^{i\theta}z))\\
&=&z+\frac{K(z)}{K'(z)}t+iz\theta +(\frac{K(z)}{K'(z)}-\frac{1}{2}\frac{K(z)^{2}K''(z)}{(K'(z))^{3}})t^{2}\\
&-&\frac{1}{2}z\theta^{2}+\frac{1}{2}(1-\frac{K(z)}{K'(z)}\frac{K''(z)}{K'(z)})izt\theta+O(t^{3})+O(\theta^{3})\\
&=&z+\frac{z(1+z)}{(1-z)}t+iz\theta+\frac{1}{2}\frac{z(1+z)(1+2z-z^{2})}{(1-z)^{3}}t^{2}\\&-&\frac{1}{2}z\theta^{2}
+\frac{1}{2}\frac{1+2z-z^2}{(1-z)^2}izt\theta+O(t^3)+O(\theta^3).
\end{eqnarray*}
\par On the other hand, as it is well known, the Bergman kernel of the unit
disc $\mathbb{D}$ on the diagonal is
$$K_{\mathbb{D}}(w,w)=\frac{1}{\pi}
\frac{1}{(1-|w|^{2})^{2}}.$$ Then, by Proposition 2.1, the Bergman
kernels of $D_{\zeta}$ on the diagonal are given by
$$K_{D_{\zeta}}(z,z)=\frac{1}{\pi}
\frac{1}{(1-|E_{\zeta}^{-1}(z)|^{2})^{2}}|(E_{\zeta}^{-1}(z))_{z}|^{2}.$$
Furthermore, we can calculate that
\begin{eqnarray*}
&\displaystyle{\lim_{\zeta\rightarrow 1}}&
\frac{\partial^2\log(1-|E_{\zeta}^{-1}(z)|^{2})}{\partial \zeta
 \partial\overline{\zeta}}=\frac{-2|z|^2\Re\left((1+z)/(1-z)\right)}{(1-|z|^2)},\\
&\displaystyle{\lim_{\zeta\rightarrow
1}}&\frac{\partial^2\log|(E_{\zeta}^{-1}(z))_{z}|^2}{\partial
\zeta
 \partial\overline{\zeta}}=\frac{1}{4}\Re\left(\frac{1-3z^2}{(1-z)^2}\right).
\end{eqnarray*}
Now let $z=re^{i\theta}\in D_{\zeta}$, we conclude that
$$\lim_{\zeta\rightarrow 1}\frac{\partial ^{2}\log K_{D_{\zeta}}(z,z)}{\partial
\zeta
\partial \overline{\zeta}} =\frac{1}{4}\frac{1+r^2\cos 2\theta}{1+r^2-2r\cos \theta}. $$
Then
$$\lim_{r\rightarrow 1}\lim_{\zeta\rightarrow 1}\frac{\partial ^{2}\log K_{D_{\zeta}}(z,z)}{\partial
\zeta
\partial \overline{\zeta}} =\frac{1}{4}\left((1-\cos \theta)+\frac{1}{(1-\cos \theta)}-2\right).\eqno(3.3.1)$$
If $\zeta \rightarrow 1$ then the singular point of
the boundary $D_{\zeta}$ tends to $z=1$. We see from (3.3.1) that
\newline (1) $\partial^2\log K_{D_{\zeta}}(z,z)/\partial \zeta\partial
\overline{\zeta}$ tends to $\infty$ with order 2 as
$(\zeta,z)\in\mathcal{D}$ tends to $(1,1)$,
\newline (2) $\partial^2\log K_{D_{\zeta}}(z,z)/\partial \zeta\partial
\overline{\zeta}$ tends to $0$ with order 2 as
$(\zeta,z)\in\mathcal{D}$ tends to

$(1,\pm i)$,
\newline (3) $\partial^2\log K_{D_{\zeta}}(z,z)/\partial \zeta\partial
\overline{\zeta}$ tends to a positive number as
$(\zeta,z)\in\mathcal{D}$ tends to

$(1,z)\in \partial\mathcal{D}$ here $z\neq 1,\pm i$.
\end{proof}
\subsection{Case of The Family of Rectangles}
In this subsection, Theorem 1.4 and 1.5 which are the explicit expressions of Bergman kernels of rectangles and the
boundary behavior of the Levi form of $\log K_{R_\zeta}(\zeta,z)$ with respect to $\zeta$, will be proved.
\begin{proof}[Proof of Theorem 1.4] Firstly, for symmetry, we consider the
transformation $F(w,\zeta)$ which maps the upper half of the $w$
plane $\mathbb{H}$ onto $R'_{\zeta}$ which is a rectangle with
vertices $A_1(\Re\zeta)$, $A_2(\zeta)$, $A_3(-\overline{\zeta})$,
$A_4(-\Re\zeta)$ for each $\zeta$
 in the $z$ plane.
We associate $A_1(\Re\zeta)$ with $a_1(1)$, $A_2(\zeta)$ with
$a_2(1/k(\zeta))$, and $w=0$ with $z=0$. Then by symmetry,
$A_3(-\overline{\zeta})$, $A_4(-\Re\zeta)$ are associated with
$a_3(-1/k(\zeta))$, $a_4(-1)$ respectively. Our goal is to
determine both the transformation $z=F(w,\zeta)$ and the constant
$k$ as an analytic function with respect to the variable $\zeta$.
In this case $\alpha_1=\alpha_2=\alpha_3=\alpha_4=1/2$, $a_1=1$,
$a_2=1/k$, $a_3=-1/k$, $a_4=-1$. Furthermore, because
$F(0,\zeta)=0$ (symmetry), the constant $C'$ of integration
(2.2.1) is zero; thus (2.2.1) yields
\begin{eqnarray}
z=F(w,\zeta)=C(\zeta)\int_0^w((1-t^2)(1-k^2(\zeta)t^2))^{-\frac{1}{2}}dt.
\end{eqnarray}
The integral appearing in (3.4.1), with the choice of a single
branch defined by the requirement that $0<\arg(w-a_i)<\pi$,
$i=1,2,3,4$, is the so-called elliptic integral of the first kind.
The association of $A_1(\Re\zeta)$ with $a_1(1)$ and $A_2(\zeta)$
with $a_2(1/k(\zeta))$ imply that
\begin{eqnarray}
\Re\zeta&=&C(\zeta)\int_0^1((1-t^2)(1-k^2(\zeta)t^2))^{-\frac{1}{2}}dt,\\
\Im\zeta&=&C(\zeta)\int_0^1((1-t^2)(1-(1-k^2(\zeta))t^2))^{-\frac{1}{2}}dt.
\end{eqnarray}
 Since $k(\zeta)$ is a real valued analytic function with
respect to $\zeta$, then, let $\zeta=1+i+\varepsilon$,
$$k(\zeta)=k_0+2\Re((a+ib)\varepsilon)+2\Re((c+id)\varepsilon^2)+2e|\varepsilon|^2+\cdots $$
in the neighborhood of the point $\zeta=1+i$.
Then,
\begin{eqnarray*}
 \left(k^2(\zeta)\right)^n&=&\left(k^2_0\right)^n+n\cdot (k_0^2)^{n-1}4k_0\left(\Re((a+ib)\varepsilon)+\Re((c+id)\varepsilon^2)+e|\varepsilon|^2\right.\\
    &\ \ \ \ +&\left.1/k_0(\Re((a+ib)\varepsilon))^2\right)\\
    &\ \ \ \ +&\frac{n(n-1)}{2}(1-k_0^2)^{n-2}16k_0^2\left(\Re((a+ib)\varepsilon)\right)^2+\cdots,\\
  \left(1-k^2(\zeta)\right)^n&=&\left(1-k^2_0\right)^n-n\cdot (1-k_0^2)^{n-1}4k_0\left(\Re((a+ib)\varepsilon)+\Re((c+id)\varepsilon^2)
    \right.\\
    &\ \ \ \ +&\left.e|\varepsilon|^2+1/k_0(\Re((a+ib)\varepsilon))^2\right)\\
    &\ \ \ \ +&\frac{n(n-1)}{2}(1-k_0^2)^{n-2}16k_0^2(\Re((a+ib)\varepsilon))^2+\cdots.
\end{eqnarray*}
By the power series expression of $K(k)$ and (3.4.2), (3.4.3) we deduce that
\begin{eqnarray}
\Re\zeta=C(\zeta)\frac{\pi}{2}\left(1+\sum_{n=1}^\infty\frac{((2n-1)!!)^2}{((2n)!!)^2}k^{2n}\right),
\end{eqnarray}
\begin{eqnarray}
\Im\zeta=C(\zeta)\frac{\pi}{2}\left(1+\sum_{n=1}^\infty\frac{((2n-1)!!)^2}{((2n)!!)^2}\left(1-k^2\right)^n\right).
\end{eqnarray}
Then, (3.4.4) and (3.4.5) imply
$$\Re\zeta\cdot\left(1+\sum_{n=1}^\infty\frac{((2n-1)!!)^2}{((2n)!!)^2}\left(1-k^2\right)^n\right)
=\Im\zeta\cdot\left(1+\sum_{n=1}^\infty\frac{((2n-1)!!)^2}{((2n)!!)^2}k^{2n}\right) \eqno (3.4.6)$$
If $\zeta=1+i$, then, from (3.4.6) we deduce that $k_0=\sqrt{2}/2$.
\newline Let $\zeta=1+i+\varepsilon$, then (3.4.6) changes to
\begin{equation*}
\begin{split}
&(1+\Re\varepsilon)\cdot\Biggl(1+\sum_{n=1}^\infty\frac{((2n-1)!!)^2}{((2n)!!)^2}
\Bigl(\left(\frac{1}{2}\right)^n-n\cdot \left(\frac{1}{2}\right)^{n-1}2\sqrt{2}\Bigl(\Re((a+ib)\varepsilon) \\
&\qquad +\Re((c+id)\varepsilon^2)
    +e|\varepsilon|^2
    +\sqrt{2}(\Re((a+ib)\varepsilon))^2\Bigr)+\cdots\Bigr)\Biggr) \\
=&(1+\Im\varepsilon)\cdot\Biggl(1+\sum_{n=1}^\infty\frac{((2n-1)!!)^2}{((2n)!!)^2}\Bigl(\left(\frac{1}{2}\right)^n+n\cdot \left(\frac{1}{2}\right)^{n-1}2\sqrt{2}\Bigl(\Re((a+ib)\varepsilon)\\
&+\Re((c+id)\varepsilon^2)
    +e|\varepsilon|^2
    +\sqrt{2}(\Re((a+ib)\varepsilon))^2\Bigr)+\cdots\Bigr)\Biggr)
\end{split}
\end{equation*}
Comparing the coefficients of the first order of $\varepsilon$ on both sides of the above equation, we have $$a=b=\frac{1+\sum_{n=1}^\infty\frac{((2n-1)!!)^2}{((2n)!!)^2}\left(\frac{1}{2}\right)^n}
{\sum_{n=1}^\infty\frac{((2n-1)!!)^2}{((2n)!!)^2}n\left(\frac{1}{2}\right)^{n-1}}=\frac{K}{\left(4\sqrt{2}(2E-K)\right)}.$$
 Here, $K$ is the value of the complete elliptic integral
of the first kind at the point $k=1/\sqrt{2}$, and $E$ is the
value of the complete elliptic integral of the second kind at the
point $k=1/\sqrt{2}$.
 And comparing the coefficients of the second order of $\varepsilon$, we get that
 $$c=-a/2,\ \ \ \ \ d=e=-\sqrt{2}a^2.$$
From (3.4.1) we have
$C(\zeta)=\Re\zeta/K(k(\zeta)).$
 In summary, the transformation $F(w,\zeta)$ is given by
$$z=F(w,\zeta)=\frac{\Re\zeta}{K(k)}\int_0^w\left((1-t^2)\left(1-k^2(\zeta)t^2\right)\right)^{-\frac{1}{2}}dt\;.$$
\par Secondly, the inverse of the integral in $F(w,\zeta)$ gives $w$ as
a function of $z$ via one of the so-called Jacobi's elliptic
function $\sn(u,k)$. Then the inverse of $F(w,\zeta)$ with respect
to the first variable is given by
$$w=f(z,\zeta)=\sn\left(\frac{K(k)}{\Re\zeta}z,k(\zeta)\right).$$
 Moreover,
$$\overline{F(-\overline{w},\zeta)}=-F(w,\zeta),$$
 implies that $f(z,\zeta)$ maps $R_{\zeta}$ to $\{{w\in \mathbb{C}_w|w=a+ib, a>0,b>0}
 \}$, thus $f^2$ maps $R_{\zeta}$ to the upper half $w$ plane $\mathbb{H}$.
In addition, it is well-known that the Bergman kernel of the upper half
$w$ plane on the diagonal is
$$K_\mathbb{H}(w,w)=\frac{1}{4\pi(\Im w)^2}.$$
Then by Proposition 2.1, the Bergman kernels of $R_{\zeta}$ on the
diagonal are
\begin{eqnarray*}
K_{R_{\zeta}}(z,z)= \frac{1}{\pi (\Im
\sn^2(u,k))^2}|\sn^2(u,k)\cn^2(u,k)\dn^2(u,k)|\left|\frac{K(k)}{\Re\zeta}\right|^2
\end{eqnarray*} here $u=K(k)z/\Re\zeta$.
\end{proof}

 \begin{proof}[Proof of Theorem 1.5] From the expression of
$K_{R_{\zeta}}(z,z)$,
\begin{eqnarray*}
 \frac{\partial^2\log K_{R_{\zeta}}(z,z)}{\partial {\zeta}
 \partial\overline{\zeta}}=-2A+2B+2C.
\end{eqnarray*}
where
\begin{eqnarray*}
A&:=&\frac{\partial^2\log(\Im\sn^2(u,k))}{\partial
{\zeta}\partial\overline{\zeta}}=2\frac{\Im\left(\sn(u,k)\frac{\partial^2\sn(u,k)}{\partial \zeta\partial\overline{\zeta}}+\frac{\partial\sn(u,k)}{\partial\zeta}\frac{\partial\sn(u,k)}{\partial\overline{\zeta}}\right)}{\Im \sn^2(u,k)}\\
&+&\frac{2\Re\left(\sn^2(u,k)\frac{\partial\sn(u,k)}{\partial\zeta}\frac{\partial\sn(u,k)}{\partial\overline{\zeta}}\right)
-|\sn(u,k)|^2\left(|\frac{\partial\sn(u,k)}{\partial\zeta}|^2+|\frac{\partial\sn(u,k)}{\partial\overline{\zeta}}|^2\right)}{(\Im \sn^2(u,k))^2},\\
B&:=&\Re\frac{\partial^2\log(\sn(u,k)\cn(u,k)\dn(u,k))}{\partial \zeta
\partial\overline{\zeta}}\\
&=&\Re
\frac{\partial^2\log\sn(u,k)}{\partial \zeta
\partial\overline{\zeta}}+\Re
\frac{\partial^2\log\cn(u,k)}{\partial \zeta
\partial\overline{\zeta}}+\Re
\frac{\partial^2\log\dn(u,k)}{\partial \zeta
\partial\overline{\zeta}}, \\
C&:=&\frac{\partial^2(-\log(\Re{\zeta})+\log|K(k)|)}{\partial
{\zeta}
\partial\overline{\zeta}}.
\end{eqnarray*}
Using the expression of $k(\zeta)$ we get that
$$\lim_{\zeta\rightarrow 1+i}\frac{\partial k}{\partial \zeta}=(1+i)a,\
\lim_{\zeta\rightarrow 1+i}\frac{\partial k}{\partial \overline{\zeta}}=(1-i)a,$$
$$\lim_{\zeta\rightarrow 1+i}\frac{\partial k}{\partial\zeta}\frac{\partial k}{\partial\overline{\zeta}}=2a^2,\lim_{\zeta\rightarrow1+i}\frac{\partial^2k}{\partial\zeta\partial\overline{\zeta}}
=-2\sqrt{2}a^2.$$
and since $u=\sn^{-1}(1,k(\zeta))z/\Re\zeta$, then,
$$\lim_{\zeta\rightarrow 1+i}\frac{\partial u}{\partial \zeta}=(-1+i)\frac{u}{4},\ \
\lim_{\zeta\rightarrow 1+i}\frac{\partial u}{\partial \overline{\zeta}}=(-1-i)\frac{u}{4},\ \ \lim_{\zeta\rightarrow 1+i}\frac{\partial u}{\partial\zeta}\frac{\partial u}{\partial\overline{\zeta}}=\frac{1}{8}u^2,$$
$$\lim_{\zeta\rightarrow1+i}\frac{\partial^2u}{\partial\zeta\partial\overline{\zeta}}
=\left(4a^2+\frac{1}{4}\right)u,\ \ \lim_{\zeta\rightarrow 1+i}\left(\frac{\partial u}{\partial \zeta}\frac{\partial k}{\partial \overline{\zeta}}+\frac{\partial k}{\partial \zeta}\frac{\partial u}{\partial \overline{\zeta}}\right)
=0.$$
\textbf{Case 1.} We consider the boundary point $(1+i,0)$. In this case, $u$ tends to $0$ when $z$ tends to $0$.
Using the power series expansion of $\sn(u,k)$ in a neighborhood of $u=0$,
the following further results can now be verified,
\begin{eqnarray*}
&\displaystyle{\lim_{\zeta\rightarrow 1+i}}&\frac{\partial\sn(u,k)}{\partial\zeta}=(-1+i)\frac{u}{4}+O(u^3),\\
&\displaystyle{\lim_{\zeta\rightarrow 1+i}}&\frac{\partial\sn(u,k)}{\partial\overline{\zeta}}=(-1-i)\frac{u}{4}+O(u^3),\\
&\displaystyle{\lim_{\zeta\rightarrow 1+i}}&\frac{\partial \sn(u,k)}{\partial\zeta}\frac{\partial \sn(u,k)}{\partial\overline{\zeta}}=\frac{1}{8}u^2+O(u^4),\\
&\displaystyle{\lim_{\zeta\rightarrow 1+i}}&\frac{\partial^2\sn(u,k)}{\partial\zeta\partial\overline{\zeta}}=\left(4a^2+\frac{1}{4}\right)u+O(u^3),\\
&\displaystyle{\lim_{\zeta\rightarrow 1+i}}&\left(\sn(u,k)\frac{\partial^2\sn(u,k)}{\partial\zeta\partial\overline{\zeta}}+\frac{\partial \sn(u,k)}{\partial\zeta}\frac{\partial \sn(u,k)}{\partial\overline{\zeta}}\right)=\left(4a^2+\frac{1}{8}\right)u^2+O(u^4),\\
&\displaystyle{\lim_{\zeta\rightarrow 1+i}}&\left(\sn(u,k)\frac{\partial^2\sn(u,k)}{\partial\zeta\partial\overline{\zeta}}-\frac{\partial \sn(u,k)}{\partial\zeta}\frac{\partial \sn(u,k)}{\partial\overline{\zeta}}\right)=\left(4a^2+\frac{3}{8}\right)u^2+O(u^4),\\
&\displaystyle{\lim_{\zeta\rightarrow 1+i}}&\left(\sn^2(u,k)\frac{\partial \sn(u,k)}{\partial\zeta}\frac{\partial \sn(u,k)}{\partial\overline{\zeta}}\right)=\frac{1}{8}u^4+O(u^6),\\
&\displaystyle{\lim_{\zeta\rightarrow 1+i}}&\left(|\frac{\partial\sn(u,k)}{\partial\zeta}|^2+|\frac{\partial\sn(u,k)}{\partial\overline{\zeta}}|^2\right)
=\frac{1}{4}|u|^2+O(u^3).
\end{eqnarray*}
Also, by using the power series expansions of $\cn(u,k)$ and $\dn(u,k)$
we obtain
\begin{eqnarray*}
&\displaystyle{\lim_{\zeta\rightarrow 1+i}}&\left(\cn(u,k)\frac{\partial^2\cn(u,k)}{\partial\zeta\partial\overline{\zeta}}-\frac{\partial \cn(u,k)}{\partial\zeta}\frac{\partial \cn(u,k)}{\partial\overline{\zeta}}\right)=-\left(4a^2+\frac{3}{8}\right)u^2+O(u^4),\\
&\displaystyle{\lim_{\zeta\rightarrow 1+i}}&\left(\dn(u,k)\frac{\partial^2\dn(u,k)}{\partial\zeta\partial\overline{\zeta}}-\frac{\partial \dn(u,k)}{\partial\zeta}\frac{\partial \dn(u,k)}{\partial\overline{\zeta}}\right)=-\left(2a^2+\frac{3}{16}\right)u^2+O(u^4).
\end{eqnarray*}
It follows that
\begin{equation*}
\lim_{z\rightarrow 0,{\zeta}\rightarrow
1+i}A=8a^2+\frac{1}{4},\ \
\lim_{z\rightarrow 0,{\zeta}\rightarrow
1+i}B=4a^2+\frac{1}{8},\ \
\lim_{z\rightarrow 0,{\zeta}\rightarrow
1+i}C=\frac{1}{8}+4a^2.
\end{equation*}
Finally, we get
$$\lim_{z\rightarrow 0,{\zeta}\rightarrow 1+i}\frac{\partial^2\log K_{R_{\zeta}}(z,z)}{\partial
\zeta
\partial\overline{\zeta}}=\lim_{z\rightarrow 0,{\zeta}\rightarrow
1+i}(-2A+2B+2C)=0.$$
\textbf{Case 2.} We consider boundary point $(1+i,1)$.
By the proof of Theorem 1.3, the periods of $\sn(u,k)$ where $u=K(k)z/\Re\zeta$, are then seen to be $4\Re\zeta$ and $2i\Im\zeta$.
In such a case, $u$ tends to $K(k)$ when $z$ tends to $\Re\zeta$. The following identity
 $$f(u',k):=\sn(u'+K(k))=\frac{\cn u'}{\dn u'},$$
 can be verified by Addition Theorem for $\sn u$ and the definition of $\dn u$.
$u'$ tends to 0 when $u$ tends to $K(k)$. For simplicity, we still use $u$ to replace $u'$.
The power series expansion of $f(u,k)$ in a neighborhood of $u=0$ is
$$f(u,k)=1+\frac{1}{2}(-1+k^2)u^2+\frac{1}{24}(1-6k^2+5k^4)u^4+O(u^6).$$
Then,
\begin{eqnarray*}
&\displaystyle{\lim_{\zeta\rightarrow 1+i}}&\frac{\partial f(u,k)}{\partial\zeta}=\frac{1}{2}(1-i)\frac{u^2}{4}+\frac{\sqrt{2}}{2}(1+i)au^2+O(u^4),\\
&\displaystyle{\lim_{\zeta\rightarrow 1+i}}&\frac{\partial f(u,k)}{\partial\overline{\zeta}}=\frac{1}{2}(1+i)\frac{u^2}{4}+\frac{\sqrt{2}}{2}(1-i)au^2+O(u^4),\\
&\displaystyle{\lim_{\zeta\rightarrow 1+i}}&\frac{\partial f(u,k)}{\partial\overline{\zeta}}
\frac{\partial f(u,k)}{\partial\zeta}=\left(a^2+\frac{1}{32}\right)u^4+O(u^6),\\
&\displaystyle{\lim_{\zeta\rightarrow 1+i}}&\frac{\partial^2 f(u,k)}{\partial\zeta\partial\overline{\zeta}}
=\left(-2a^2-\frac{3}{8}\right)u^2+O(u^4),
\end{eqnarray*}
\begin{eqnarray*}
&\displaystyle{\lim_{\zeta\rightarrow 1+i}}&\left(f(u,k)\frac{\partial^2f(u,k)}{\partial\zeta\partial\overline{\zeta}}+\frac{\partial f(u,k)}{\partial\zeta}\frac{\partial f(u,k)}{\partial\overline{\zeta}}\right)=\left(-2a^2-\frac{3}{8}\right)u^2+O(u^4),\\
&\displaystyle{\lim_{\zeta\rightarrow 1+i}}&\left(f^2(u,k)\frac{\partial f(u,k)}{\partial\zeta}\frac{\partial f(u,k)}{\partial\overline{\zeta}}\right)=\left(a^2+\frac{1}{32}\right)u^4+O(u^6),\\
&\displaystyle{\lim_{\zeta\rightarrow 1+i}}&\left(|\frac{\partial f(u,k)}{\partial\zeta}|^2+|\frac{\partial f(u,k)}{\partial\overline{\zeta}}|^2\right)
=\left(2a^2+\frac{1}{16}\right)|u|^4+O(u^6).
\end{eqnarray*}
The following identities can be established from (2.4.1)--(2.4.7),
$$ \cn(u+K(k))=-k'\frac{\sn u}{\dn u},\ \ \ \dn(u+K(k))=k'\frac{1}{\dn u}.$$
These identities and expressions of $A,B,C$ give
\begin{eqnarray*}
A&:=&\frac{\partial^2\log(\Im\sn^2(u+K(k),k))}{\partial
{\zeta}\partial\overline{\zeta}}=2\frac{\Im\left(f(u,k)\frac{\partial^2f(u,k)}{\partial \zeta\partial\overline{\zeta}}+\frac{\partial f(u,k)}{\partial\zeta}\frac{\partial f(u,k)}{\partial\overline{\zeta}}\right)}{\Im f^2(u,k)}\\
&+&\frac{2\Re\left(f^2(u,k)\frac{\partial f(u,k)}{\partial\zeta}\frac{\partial f(u,k)}{\partial\overline{\zeta}}\right)
-|f(u,k)|^2\left(|\frac{\partial f (u,k)}{\partial\zeta}|^2+|\frac{\partial f(u,k)}{\partial\overline{\zeta}}|^2\right)}{(\Im f^2(u,k))^2},\\
B&:=&\Re\frac{\partial^2\log(\sn(u+K(k),k)\cn(u+K(k),k)\dn(u+K(k),k))}{\partial \zeta
\partial\overline{\zeta}}\\
&=&\log(k'^2)+\Re
\frac{\partial^2\log\sn(u,k)}{\partial \zeta
\partial\overline{\zeta}}+\Re
\frac{\partial^2\log\cn(u,k)}{\partial \zeta
\partial\overline{\zeta}}-3\Re
\frac{\partial^2\log\dn(u,k)}{\partial \zeta
\partial\overline{\zeta}}, \\
C&:=&\frac{\partial^2(-\log(\Re{\zeta})+\log|K(k)|)}{\partial
{\zeta}
\partial\overline{\zeta}}.
\end{eqnarray*}
These immediately induce that,
\begin{equation*}
\lim_{z\rightarrow 1,{\zeta}\rightarrow
1+i}A=-8a^2-\frac{1}{2},\ \
\lim_{z\rightarrow 1,{\zeta}\rightarrow
1+i}B=\frac{1}{8}-12a^2,\ \
\lim_{z\rightarrow 1,{\zeta}\rightarrow
1+i}C=\frac{1}{8}+4a^2.
\end{equation*}
Hence, the following result emerges,
$$\lim_{z\rightarrow 1}\left(\lim_{{\zeta}\rightarrow 1+i}\frac{\partial^2\log K_{R_{\zeta}}(z,z)}{\partial
\zeta
\partial\overline{\zeta}}\right)=\lim_{z\rightarrow 1,{\zeta}\rightarrow
1+i}(-2A+2B+2C)=\frac{3}{2}.$$
\textbf{Case 3.} We consider the boundary point $(1+i,i)$. In this circumstance, $u$ tends to $iK'(k)$ when $z$ tends to $i\Im\zeta$ here $K'(k)=K(k')$ defined by
$$K'(k)=\int_0^1((1-t^2)(1-(1-k^2(\zeta))t^2))^{-\frac{1}{2}}dt.$$
The following identities can be established from (2.4.1)-(2.4.14),
$$\sn(u+iK'(k))=\frac{1}{k}\frac{1}{\sn u},\;\;\;\;
\cn(u+iK'(k))=\frac{1}{ik}\frac{\dn u}{\sn
u},\;\;\;\;\dn(u+iK'(k))=-i\frac{\cn u}{\sn u}.$$
The power series expansion of $f(u,k):=1/(k\sn(u,k))$ is
$$f(u,k)\equiv\frac{1}{k}\frac{1}{\sn(u,k)}=\frac{1}{ku}+\frac{1}{6}\left(k+\frac{1}{k}\right)u
          +\frac{1}{360}\left(7k^3-22k+\frac{7}{k}\right)u^3+O(u^5).$$
 Then, by a similar process as in Case 2:
 \begin{eqnarray*}
&\displaystyle{\lim_{\zeta\rightarrow 1+i}}&\frac{\partial f(u,k)}{\partial\zeta}=\left(\frac{\sqrt{2}}{4}(1-i)-2(1+i)a\right)\frac{1}{u}+O(u),\\
&\displaystyle{\lim_{\zeta\rightarrow 1+i}}&\frac{\partial f(u,k)}{\partial\overline{\zeta}}=\left(\frac{\sqrt{2}}{4}(1+i)-2(1-i)a\right)\frac{1}{u}+O(u),\\
&\displaystyle{\lim_{\zeta\rightarrow 1+i}}&\frac{\partial f(u,k)}{\partial\overline{\zeta}}
\frac{\partial f(u,k)}{\partial\zeta}=\left(\frac{1}{4}+8a^2\right)\frac{1}{u^2}+O(1),\\
&\displaystyle{\lim_{\zeta\rightarrow 1+i}}&\frac{\partial^2 f(u,k)}{\partial\zeta\partial\overline{\zeta}}
= 4\sqrt{2}a^2\cdot\frac{1}{u}+O(u),\\
&\displaystyle{\lim_{\zeta\rightarrow 1+i}}&\left(f(u,k)\frac{\partial^2f(u,k)}{\partial\zeta\partial\overline{\zeta}}+\frac{\partial f(u,k)}{\partial\zeta}\frac{\partial f(u,k)}{\partial\overline{\zeta}}\right)=\left(\frac{1}{4}+16a^2\right)\frac{1}{u^2}+O(1),\\
&\displaystyle{\lim_{\zeta\rightarrow 1+i}}&\left(f^2(u,k)\frac{\partial f(u,k)}{\partial\zeta}\frac{\partial f(u,k)}{\partial\overline{\zeta}}\right)=\left(\frac{1}{2}+16a^2\right)\frac{1}{u^2}+O(1),\\
&\displaystyle{\lim_{\zeta\rightarrow 1+i}}&\left(|\frac{\partial f(u,k)}{\partial\zeta}|^2+|\frac{\partial f(u,k)}{\partial\overline{\zeta}}|^2\right)
=\left(\frac{1}{2}+16a^2\right)\frac{1}{u^2}+O(1).
\end{eqnarray*}
From these results we may verify that the result listed below follows:
\begin{equation*}
\lim_{z\rightarrow i,{\zeta}\rightarrow
1+i}A=-\frac{1}{4},\ \
\lim_{z\rightarrow i,{\zeta}\rightarrow
1+i}B=4a^2-\frac{3}{8},\ \
\lim_{z\rightarrow i,{\zeta}\rightarrow
1+i}C=\frac{1}{8}+4a^2.
\end{equation*}
Then,
$$\lim_{z\rightarrow i}\left(\lim_{{\zeta}\rightarrow 1+i}\frac{\partial^2\log K_{R_{\zeta}}(z,z)}{\partial
\zeta
\partial\overline{\zeta}}\right)=\lim_{z\rightarrow i,{\zeta}\rightarrow
1+i}(-2A+2B+2C)=16a^2>0.$$
\textbf{Case 4.} We consider the boundary point $(1+i,1+i)$. In this instance, $u$ tends to $K(k)+iK'(k)$ when $z$ tends to $\zeta$.
The following identities can also be established from (2.4.1)-(2.4.7),
\begin{eqnarray*}
\sn(u+K(k)+iK'(k))&=&\frac{1}{k}\frac{\dn u}{\cn u},\ \ \ \ \
\cn(u+K(k)+iK'(k))=\frac{k'}{ik}\frac{1}{\cn
u},\\
\dn(u+K(k)+iK'(k))&=&ik'\frac{\sn u}{\cn u}.
\end{eqnarray*}
The power series expansion of $f(u,k):=\dn (u,k)/(k\cn(u,k))$ is
$$f(u,k)\equiv\frac{1}{k}\frac{\dn u}{\cn u}=\frac{1}{k}+\frac{1}{2}\left(-k+\frac{1}{k}\right)u^2
          +\frac{1}{24}\left(k^3-6k+\frac{5}{k}\right)u^4+O(u^6).$$
 Then, applying the same procedure as in Case 2, we get that
 \begin{eqnarray*}
&\displaystyle{\lim_{\zeta\rightarrow 1+i}}&\frac{\partial f(u,k)}{\partial\zeta}=(-2)(1+i)a+\left(\frac{\sqrt{2}}{8}(-1+i)-\frac{3}{2}(1+i)a\right)u^2+O(u^4),\\
&\displaystyle{\lim_{\zeta\rightarrow 1+i}}&\frac{\partial f(u,k)}{\partial\overline{\zeta}}=(-2)(1-i)a+\left(\frac{\sqrt{2}}{8}(-1-i)-\frac{3}{2}(1-i)a\right)u^2+O(u^4),\\
&\displaystyle{\lim_{\zeta\rightarrow 1+i}}&\frac{\partial f(u,k)}{\partial\overline{\zeta}}
\frac{\partial f(u,k)}{\partial\zeta}=8a^2+12a^2\cdot u^2+O(u^2),\\
&\displaystyle{\lim_{\zeta\rightarrow 1+i}}&\frac{\partial^2 f(u,k)}{\partial\zeta\partial\overline{\zeta}}
= 12\sqrt{2}a^2+\left(\frac{3\sqrt{2}}{16}+9\sqrt{2}a^2\right)\cdot u^2+O(u^4),\\
&\displaystyle{\lim_{\zeta\rightarrow 1+i}}&\left(f(u,k)\frac{\partial^2f(u,k)}{\partial\zeta\partial\overline{\zeta}}+\frac{\partial f(u,k)}{\partial\zeta}\frac{\partial f(u,k)}{\partial\overline{\zeta}}\right)\\
&=&32a^2+\left(\frac{3}{8}+36a^2\right)u^2+O(u^4),\\
&\displaystyle{\lim_{\zeta\rightarrow 1+i}}&\left(f^2(u,k)\frac{\partial f(u,k)}{\partial\zeta}\frac{\partial f(u,k)}{\partial\overline{\zeta}}\right)=32a^2+64a^2\cdot u^2+O(u^4),\\
&\displaystyle{\lim_{\zeta\rightarrow 1+i}}&|f(u,k)|^2=\displaystyle{\lim_{\zeta\rightarrow 1+i}} f(u,k)\overline{f(u,k)}=2+\frac{\sqrt{2}}{2}\Re u^2+O(u^4),\\
&\displaystyle{\lim_{\zeta\rightarrow 1+i}}&\left(|\frac{\partial f(u,k)}{\partial\zeta}|^2+|\frac{\partial f(u,k)}{\partial\overline{\zeta}}|^2\right)
=\displaystyle{\lim_{\zeta\rightarrow 1+i}}2\left(\frac{\partial f(u,k)}{\partial\zeta}\frac{\partial \overline{f(u,k)}}{\partial\overline{\zeta}}\right)\\
&=&8a^2+12a^2\Re u^2-\sqrt{2}\Im u^2+O(u^4).
\end{eqnarray*}
Then,
\begin{eqnarray*}
&\displaystyle{\lim_{z\rightarrow 1+i,{\zeta}\rightarrow
1+i}}&A=\displaystyle{\lim_{u\rightarrow 0}}\left(-72a^2-\frac{3}{4}+\frac{4\sqrt{2}a^2}{\Im u^2}+\left(16-8\sqrt{2}\right)a^2\frac{\Re u^2}{(\Im u^2)^2}+O(1)\right),\\
&\displaystyle{\lim_{z\rightarrow 1+i,{\zeta}\rightarrow
1+i}}&B=4a^2+\frac{1}{8},\ \ \ \ \
\displaystyle{\lim_{z\rightarrow 1+i,{\zeta}\rightarrow
1+i}}C=4a^2+\frac{1}{8}.
\end{eqnarray*}
Finally,
$$\lim_{z\rightarrow 1+i}\left(\lim_{{\zeta}\rightarrow 1+i}\frac{\partial^2\log K_{R_{\zeta}}(z,z)}{\partial
\zeta
\partial\overline{\zeta}}\right)=\lim_{z\rightarrow 1+i,{\zeta}\rightarrow
1+i}(-2A+2B+2C)=\infty.$$

\end{proof}
\subsection{Case of The Family of Half Strips}

\begin{proof}[Proof of Theorem 1.6]
 Applying the same strategy as Theorem
1.3, the biholomorphic mappings between the upper half $w$ plane
and $S_\zeta$ are
$$w=f^2(z,\zeta)=\sin^2 u,$$ where $u=\pi z/2\Re \zeta$. Here we associate $A_1(0)$ with $a_1(0)$, $A_2(\Re\zeta)$ with
$a_2(1)$, and $A(\infty)$ with $a(\infty)$. By Proposition 2.1,
the Bergman kernels of $S_\zeta$ are
$$K_{S_\zeta}(z,z)=
\frac{1}{\pi(\Im \sin^2u)^2}\left|\frac{\pi}{2\Re\zeta}\sin u\cos
u \right|^2.$$
Since $u=\pi z/2\Re \zeta$, then
$$\lim_{\zeta\rightarrow 1}\frac{\partial u}{\partial\zeta}=\lim_{\zeta\rightarrow 1}\frac{\partial u}{\partial\overline{\zeta}}=-\frac{u}{2},\ \ \; \lim_{\zeta\rightarrow 1} \frac{\partial u}{\partial\zeta}\frac{\partial u}{\partial\overline{\zeta}}=\frac{u^2}{4},\ \ \ \;\lim_{\zeta\rightarrow 1}\frac{\partial^2 u}{\partial\zeta\partial\overline{\zeta}}=\frac{u}{2}.$$
Also,
$$\frac{\partial^2\log \Im \sin^2 u}{\partial\zeta\partial\overline{\zeta}}=
\frac{\Im\left(\frac{\partial^2u}{\partial\zeta\partial\overline{\zeta}}\sin2u+2\cos2u\frac{\partial
u}{\partial\zeta}\frac{\partial
u}{\partial\overline{\zeta}}\right)}{\Im \sin^2u}-\left(\frac{\Im
\left(\frac{\partial u}{\partial\zeta}\sin2u\right)}{\Im
\sin^2u}\right)^2, $$
using the power series expansions of $\sin u$ and $\cos u$ in a neighborhood of $u=0$ and the periodicity of $\sin u$ and $\cos u$, we have
$$\lim_{z\rightarrow 0}\left(\lim_{\zeta\rightarrow
1}\frac{\partial^2\log \Im
\sin^2u}{\partial\zeta\partial\overline{\zeta}}\right)=\frac{1}{2},\ \ \lim_{z\rightarrow
1}\left(\lim_{\zeta\rightarrow 1}\frac{\partial^2\log \Im
\sin^2u}{\partial\zeta\partial\overline{\zeta}}\right)=\frac{1}{4}.$$
Moreover,
$$\frac{\partial^2\log|\frac{\pi}{2\Re\zeta}\sin u\cos u|}{\partial\zeta\partial\overline{\zeta}}
=\Re\left(\frac{\partial^2u}{\partial\zeta\partial\overline{\zeta}}\frac{\sin
4u}{\sin^22u}-\frac{\partial u}{\partial\zeta}\frac{\partial
u}{\partial\overline{\zeta}}\frac{4}{\sin^22u}\right)+\frac{1}{(\zeta+\overline{\zeta})^2}$$
then
\begin{eqnarray*}
\displaystyle{\lim_{z\rightarrow 0}\left(\lim_{\zeta\rightarrow
1}\frac{\partial^2\log|\frac{\pi}{2\Re\zeta}\sin u\cos
u|}{\partial\zeta\partial\overline{\zeta}}\right)}&=&\frac{1}{2},\\
\displaystyle{\lim_{z\rightarrow
1}\left(\lim_{\zeta\rightarrow
1}\frac{\partial^2\log|\frac{\pi}{2\Re\zeta}\sin u\cos
u|}{\partial\zeta\partial\overline{\zeta}}\right)}&=&\frac{1}{4}\lim_{z\rightarrow
1}\frac{1}{\Re(1-z)^2}=\infty.
\end{eqnarray*}
 Hence,
 $$\lim_{z\rightarrow 0}\left(\lim_{\zeta\rightarrow
1}\frac{\partial^2\log K_{S_\zeta}(z,z)}{\partial\zeta\partial\overline{\zeta}}\right)=0,\ \ \lim_{z\rightarrow 1}\left(\lim_{\zeta\rightarrow
1}\frac{\partial^2\log K_{S_\zeta}(z,z)}{\partial\zeta\partial\overline{\zeta}}\right)=\infty.
$$
Next, we consider the boundary point $(\zeta,\infty)$. $u$ tends to $\pi z/2$ when $\zeta$ tends to 1. Let $\frac{\pi}{2}z=x+iy$ with $x\in (0,\frac{\pi}{2})$, $y>0$.
\begin{eqnarray*}
\displaystyle{\lim_{\zeta\rightarrow 1}}\,\sin u&=&\displaystyle{\lim_{\zeta\rightarrow 1}}\frac{e^{iu}-e^{-iu}}{2i}=\frac{\sin x}{2}\left(e^y+e^{-y}\right)-\frac{i}{2}\cos x\left(-e^y+e^{-y}\right),\\
\displaystyle{\lim_{\zeta\rightarrow 1}}\,\cos u&=&\displaystyle{\lim_{\zeta\rightarrow 1}}\frac{e^{iu}+e^{-iu}}{2}=\frac{\cos x}{2}\left(e^y+e^{-y}\right)+\frac{i}{2}\sin x\left(-e^y+e^{-y}\right).
\end{eqnarray*}
Then
for fixed $x$
\begin{eqnarray*}
\displaystyle{\lim_{y\rightarrow \infty}}\left(\lim_{\Re\zeta\rightarrow
1}\frac{\partial^2\log \Im \sin^2 u}{\partial\zeta\partial\overline{\zeta}}\right)&=&y-x^2+x\cot2x-x^2\cot^2x,\\
\displaystyle{\lim_{y\rightarrow \infty}}\left(\lim_{\Re\zeta\rightarrow
1}\frac{\partial^2\log|\frac{\pi}{2\Re\zeta}\sin u\cos
u|}{\partial\zeta\partial\overline{\zeta}}\right)&=&y+\frac{1}{4}.
\end{eqnarray*}
So,
$$\lim_{y\rightarrow \infty}\left(\lim_{\zeta\rightarrow
1}\frac{\partial^2\log K_{S_\zeta}(z,z)}{\partial\zeta\partial\overline{\zeta}}\right)=2x^2+2\left(x\cot
2x-\frac{1}{2}\right)^2.$$ This implies that $\lim_{y\rightarrow
\infty}\left(\lim_{\zeta\rightarrow
1}\frac{\partial^2\log K_{S_\zeta}(z,z)}{\partial\zeta\partial\overline{\zeta}}\right)$
exists for fixed $x\in (0,\frac{\pi}{2})$ and
$$\lim_{x\rightarrow 0}\left(\lim_{y\rightarrow \infty}\left(\lim_{\zeta\rightarrow
1}\frac{\partial^2\log K_{S_\zeta}(z,z)}{\partial\zeta\partial\overline{\zeta}}\right)\right)=0,$$
$$\lim_{x\rightarrow \pi/2}\left(\lim_{y\rightarrow \infty}\left(\lim_{\zeta\rightarrow
1}\frac{\partial^2\log K_{S_\zeta}(z,z)}{\partial\zeta\partial\overline{\zeta}}\right)\right)=\infty.$$
\end{proof}
We are grateful for Professor Takeo Ohsawa and my friend Masanori Adachi to their advice.
 \vskip 4mm {

}

\end{document}